\RequirePackage{fix-cm}
\documentclass[smallextended,numbook,envcountsame]{svjour3}       % onecolumn (second format)
\smartqed  % flush right qed marks, e.g. at end of proof
\usepackage[T1]{fontenc}
\usepackage{lmodern}

\usepackage{latexsym}
\usepackage{amsfonts}
\usepackage{amsmath}
\usepackage{amssymb}
\usepackage{graphicx}
\usepackage{color}
\usepackage{enumerate}

\usepackage{mathtools}

\usepackage{graphics}
\usepackage{epsfig}

\usepackage{calrsfs}
\usepackage{esint}

\usepackage%[final]
[%notcite,
final
]{showkeys}%

\usepackage{hyperref}% 

\makeatletter
\renewcommand\theequation{%\thesection.
\@arabic\c@equation}%
\makeatother

\renewenvironment{proof}[1][]{\medskip \noindent\emph{Proof}{#1} } {\qed\medskip}

\makeatletter
  \renewcommand\section{\@startsection {section}{1}{\z@}%
                                   {-\bigskipamount}%
                                   {\medskipamount}%
                                   {\large\bfseries%\mathversion{bold}
                                   \raggedright}}

  \renewcommand\subsection{\@startsection {subsection}{2}{\z@}%
                                   {-\medskipamount}%
                                   {\smallskipamount}%
                                   {\bfseries%\mathversion{bold}
                                   \raggedright}}
\makeatother

\renewcommand{\gg}{>\kern-2pt>}
\renewcommand{\ll}{<\kern-2pt<}

\renewcommand{\gg}{>\kern-2pt>}
\renewcommand{\ll}{<\kern-2pt<}

\renewcommand{\le}{\leqslant}
\renewcommand{\ge}{\geqslant}

\newcommand{\al}{\alpha}

\newcommand{\si}{\sigma}
\newcommand{\Si}{\Sigma}

\newcommand{\PP}{\operatorname{\mathcal{P}}}
\newcommand{\T}{\operatorname{\mathcal{T}}}

\newcommand{\LD}{\mathcal{L}\!\mathcal{D}}
\renewcommand{\LD}{\mathcal{L}{\kern -1.9pt}\mathcal{D}}
\renewcommand{\LD}{\mathcal{D}}
\renewcommand{\LD}{\mathcal{L}{\kern -4pt}\mathcal{C}}
\renewcommand{\LD}{\mathcal{R}{\kern -3pt}\mathcal{C}}

\newcommand{\ii}[1]{\mathrm{I}\!\left\{#1\right\}}

\newcommand{\Z}{\mathbb{Z}}
\newcommand{\R}{{\mathbb{R}}}

\renewcommand{\PP}[1]{\mathcal{P}^{#1}_+}
\renewcommand{\PP}{\mathcal{P}}

\newcommand{\p}{q}
\newcommand{\pp}{\mathbf{q}}

\def\subclassname{{\bfseries Mathematics Subject Classification
(2010)}\enspace}

\begin{document}

\title{Nonnegative sum-symmetric matrices, optimal-score partitions, and optimal resource allocation %\thanks{Grants or other notes
%about the article that should go on the front page should be
%placed here. General acknowledgments should be placed at the end of the article.}
}
%\subtitle{Do you have a subtitle?\\ If so, write it here}

\titlerunning{Optimal-score partitions}        % if too long for running head

\author{Iosif Pinelis       
% \and
%        Second Author %etc.
}

%\authorrunning{Short form of author list} % if too long for running head

\institute{I. Pinelis \at
              Department of Mathematical Sciences\\
Michigan Technological University\\
Hough\-ton, Michigan 49931, USA\\
              Tel.: +1-906-487-2108\\
              Fax: +1-906-487-3133\\
              \email{ipinelis@mtu.edu}           %  \\
%             \emph{Present address:} of F. Author  %  if needed
%           \and
%           S. Author \at
%              second address
}

\date{Received: date / Accepted: date}
% The correct dates will be entered by the editor

\maketitle

\begin{abstract}
The main result of the note describes certain optimal-score partitions, which can be interpreted as optimal resource allocations. This result is based on the fact that any nonnegative square matrix whose column sums are the same as the corresponding row sums can be represented as the sum of circuit matrices. 
% \PACS{PACS code1 \and PACS code2 \and more}
\end{abstract}

%\tableofcontents

\section{Nonnegative sum-symmetric matrices}
A matrix is called nonnegative if all its entries are nonnegative.
A square matrix $T=(t_{ij})_{i,j\in[n]}$, where $[n]:=\{1,\dots,n\}$, is called sum-symmetric if for each $i\in[n]$ the row sum $%IP04-25-19 r
s_i(T):=\sum_{j\in[n]} t_{ij}$ is the same as the corresponding column sum $c_i(T):=\sum_{j\in[n]} t_{ji}$. A matrix $(c_{ij})_{i,j\in[n]}$ is called a circuit matrix if for some set $J\subseteq[n]$, some cyclic permutation $\pi$ of $J$, and all $i,j$ in $[n]$ we have $c_{ij}=\ii{j=\pi(i)\in J}$, where $\ii\cdot$ denotes the indicator. 
Clearly, any circuit matrix is sum-symmetric. 
A central result here is that any nonnegative sum-symmetric real matrix is a conical combination of circuit matrices; see e.g.\ \cite[Theorem~1]{dantzig85} or \cite[Lemma~3.4.3]{bapat-ragh}; in \cite{dantzig85}, the sum-symmetric and circuit matrices are referred to as line-sum-symmetric and simple circuit matrices, respectively. 

This result has a short and simple proof, which extends almost verbatim to the case when the entries of the matrix are from a linearly ordered Abelian group $(G,+,0,\ge)$, with a linear order $\ge$ on the set $G$ such that for any $a$ and $b$ in $G$ one has $a\ge b\iff a-b\ge0$. Write $a>b$ to mean that $a\ge b\ne a$. It is shown in \cite{levi42} that an Abelian group can be linearly ordered iff it is torsion free, that is, iff $0$ is its only element of finite order. It is also known (see e.g.\ \cite{hahn07,gravett}) that any linearly ordered Abelian group can be embedded into the additive group $\R^I$ endowed with a lexicographical order, where $I$ is a certain linearly ordered set and $\R^I$ is the set of all functions from $I$ to $\R$ vanishing outside a well-ordered subset of $I$. Examples of linearly ordered groups are any linearly ordered rings and, in particular, any linearly ordered fields. 
%, and hence any ordered vector spaces (such spaces were apparently first introduced by Kantorovich; see e.g.\ \cite[Section ``Linear programming'']{kantorovich}. 
So, the additive groups of the ordered fields $\R$ and ${}^*\R$ of real and hyperreal numbers are linearly ordered groups. Any subgroup of any linearly ordered group is a linearly ordered group, with the inherited order. The direct product $G_1\times G_2\times\cdots$ of any linearly ordered groups $G_1,G_2,\ldots$ is a linearly ordered group with respect to the lexicographic order. 

%More specifically, let us say that an Abelian group 
%$(G,+,0)$ is \emph{nicely quasi-ordered} if $G$ is endowed with a binary relation $\ge$ %\,\subseteq G\times G$ 
%such that for all $a,b,c$ in $G$ one has (i) $a\ge b\ \&\ b\ge c\implies a\ge c$ (transitivity); (ii) $a\ge b$ or $b\ge a$ (connex property); (iii) $a\ge b\implies a-b\ge0$; and (iv) $a\ge0\ \&\ -a\ge0\implies a=0$; the reflexivity or anti-symmetry properties are not assumed here; however, it is easy to see that properties (iii) and (iv) of the relation $\ge$ imply the anti-symmetry; also, if the relation $\ge$, with properties (i)--(iv), is reflexive, then it is a linear order. 
In this group context, let us also extend the notion of a \emph{circuit matrix}, by defining it as a matrix $(c_{ij})_{i,j\in[n]}\in G^{n\times n}$ such that for some $c\in G$, some set $J\subseteq[n]$, some cyclic permutation $\pi$ of $J$, and all $i,j$ in $[n]$ we have $c_{ij}=c$ if $j=\pi(i)\in J$ and $c_{ij}=0$ otherwise; let us denote this circuit matrix by $C^{J,\pi,c}$. Now we can state  

\begin{theorem}\label{th:sum-symm}
Let $(G,+,0,\ge)$ be a linearly ordered Abelian group. 
Then any nonnegative sum-symmetric matrix in $G^{n\times n}$ is the sum of nonnegative circuit matrices in $G^{n\times n}$. 
\end{theorem}

%In the rest of the paper, we shall only need Theorem~\ref{th:sum-symm} when $G$ is $\R$ or $\Z$. 
%
%The case $G=\R$ of Theorem~\ref{th:sum-symm} complements the famous Birkhoff--von Neumann theorem, which states that every nonnegative doubly stochastic matrix is a convex combination of permutation matrices. 

For readers' convenience, let us give here 

\begin{proof}%[Proof of Theorem~\ref{th:sum-symm}]
\emph{of Theorem~\ref{th:sum-symm}.\ \ }
Take any nonnegative sum-symmetric matrix $T=(t_{ij})_{i,j\in[n]}\in G^{n\times n}$. If $s_k(T)=0$ for some $k\in[n]$, then $c_k(T)=s_k(T)=0$, and so, all entries of the $k$th row and $k$th column of $T$ are $0$. Crossing out these row and column, we obtain a nonnegative sum-symmetric matrix in $G^{(n-1)\times(n-1)}$, and the proof can be easily completed by induction on $n$. 

So, without loss of generality %(wlog) 
$s_i(T)>0$ for all $i\in[n]$, that is, for each $i\in[n]$ there is some $j\in[n]$ such that $t_{ij}>0$. Therefore, for any $i_1\in[n]$ we have a sequence $(i_1,i_2,\dots)$ in the set $[n]$ such that $t_{i_\al,i_{\al+1}}>0$ for all natural $\al$. By the pigeonhole principle, there are natural $k$ and $\ell$ with the property that $k<\ell$ and $i_k=i_\ell$. Taking such $k$ and $\ell$ with the smallest value of $\ell-k$, we will have $i_k,\dots,i_{\ell-1}$ be pairwise distinct. So, 
the condition $\pi(i_\al)=i_{\al+1}$ for $\al=k,\dots,\ell-1$ will define a cyclic permutation $\pi$ on the set $J:=\{i_k,\dots,i_{\ell-1}\}$. Then the matrix $\tilde T:=T-C^{J,\pi,t}$, where $t:=\bigwedge_{\al=k}^{\ell-1} t_{i_\al,i_{\al+1}}$, will be nonnegative and sum-symmetric, and $\tilde T$ will have strictly fewer nonzero entries than $T$ does. Now the proof can be easily completed by induction on the number of nonzero entries of the matrix. 
\end{proof}

The case $G=\R$ of Theorem~\ref{th:sum-symm} complements the famous Birkhoff--von Neumann theorem, which states that every %nonnegative 
doubly stochastic matrix is a convex combination of permutation matrices. 
One can similarly extend the Birkhoff--von Neumann theorem to groups: 

\begin{theorem}\label{th:birkhoff}
Let $(G,+,0,\ge)$ be a linearly ordered Abelian group. 
Then any nonnegative matrix $T\in G^{n\times n}$ with $s_1(T)=\cdots=s_n(T)=c_1(T)=\cdots=c_n(T)$ is the sum of nonnegative circuit matrices in $G^{n\times n}$ of the form $C^{J,\pi,c}$ with $J=[n]$. 
\end{theorem}

For a proof of Theorem~\ref{th:birkhoff}, one may take, almost verbatim (cf.\ the above proof of Theorem~\ref{th:sum-symm}), the proof of Theorem~5.1.9 in \cite{hall67}, which is based on Ph.\ Hall's theorem on distinct representatives -- see e.g.\ Theorem~5.1.1 in \cite{hall67}; other proofs of Ph.\ Hall's theorem and its extensions can be found e.g. in \cite{ann-comb} and \cite[Section~3.3]{representant}. 

In the rest of the paper, we shall only need Theorem~\ref{th:sum-symm} when $G$ is $\R$ or $\Z$.

\section{Optimal-score partitions} \label{main}
Let $k$ be a natural number. Let $\mu$ and $\nu$ be finite measures on a measurable space $(X,\Si)$ such that $\mu$ is absolutely continuous with respect to $\nu$, with a Radon--Nikodym derivative $f=\frac{d\mu}{d\nu}$. Let $\PP_k$ denote the set of all %$\Si$-measurable 
partitions $P=(A_1,\dots,A_k)$ of $X$ such that $A_i\in\Si$ for all $i\in[k]$. 

%Let $G$ be either $\R$ or $\Z$. 
Suppose that one of the following two conditions on the group $G$, the \break 
$\si$-algebra $\Si$, and the measure $\nu$ holds: 
\begin{enumerate}[(I)]
	\item \label{I} $G=\R$ and $\nu$ is non-atomic; 
	\item \label{II} $G=\Z$, $\Si$ is the powerset $2^X$ of $X$, and $\nu$ is the counting measure (so that the set $X$ is finite). 
\end{enumerate}
Then 
\begin{equation}\label{eq:p}
	\forall A\in\Si\quad\forall q\in G\cap[0,\nu(A)]\quad \exists B\in\Si\quad B\subseteq A\ \ \&\ \ \nu(B)=q.  
\end{equation}
Indeed, this is obvious when condition (II) holds. In the case when (I) holds, conclusion \eqref{eq:p} follows immediately from the well-known fact that the set of all values of a non-atomic finite measure is convex; see e.g.\ \cite[Proposition~A.1]{dudley-norv}. 

Fix any $k$-tuple %nonnegative 
\begin{equation*}
	\pp=(\p_1,\dots,\p_k)\in\big(G\cap[0,\infty)\big)^k
\end{equation*}
%$\p_1,\dots,\p_k$ in $G\cap[0,\infty)$ 
%
such that 
\begin{equation}\label{eq:sum q}
	\p_1+\dots+\p_k=\nu(X). 
\end{equation}
%in view of condition (I)$\wedge$(II), 
%such $\p_1,\dots,\p_k$ exist in the case (I) if $\nu(X)>0$, and in the case (II) if $\nu(X)\ge k$. 
Consider %now 
\begin{equation}\label{eq:P_nu}
	\PP_{\nu,\pp%\p_1,\dots,\p_k
	}:=\big\{P=(A_1,\dots,A_k)\in\PP_k\colon\nu(A_i)=\p_i\,\ \forall i\in[k]\big\}. 
\end{equation}
In view of \eqref{eq:p}, $\PP_{\nu,\pp}\ne\emptyset$. 
Moreover, let us state 

\begin{proposition}\label{prop:decr}
There exists a partition $Q=(B_1,\dots,B_k)\in\PP_{\nu,\pp}$ such that for %\break 
any $i,j$ in $[k]$ 
\begin{equation}\label{eq:inf>sup}
	\sup_{B_i}f\le\inf_{B_j}f\quad\text{whenever}\quad %i\in[k],j\in[k],
	i<j; 
\end{equation}
recall here that $\sup\emptyset=-\infty$ and $\inf\emptyset=\infty$. 
\end{proposition}

Also, fix arbitrary real numbers $s_1,\dots,s_k$ such that 
\begin{equation}\label{eq:s incr}
	s_1\le\dots\le s_k
\end{equation}
and define the ``score'' 
\begin{equation}\label{eq:s}
	s(P):=\sum_{i=1}^k s_i\mu(A_i)
\end{equation}
of any partition $P=(A_1,\dots,A_k)\in\PP_{\nu,\pp}$. 

Now we can state the main result of this note: 

\begin{theorem}\label{th:opt score}
For any partition $Q$ as in Proposition~\ref{prop:decr} and any partition $P\in\PP_{\nu,\pp}$, we have $s(Q)\ge s(P)$; that is, any partition $Q$ as in Proposition~\ref{prop:decr} has the highest possible score among all partitions in $\PP_{\nu,\pp}$. 
\end{theorem}

Let us now prove the above statements. 

\begin{proof}\emph{of Proposition~\ref{prop:decr}.\ \ } 
This will be done by induction on $k$. The case $k=1$ is trivial. By writing $\p_1+\dots+\p_k=(\p_1+\dots+\p_{k-1})+\p_k$, we reduce the consideration to the case $k=2$, so that $\pp=(\p_1, \p_2)\in\big(G\cap[0,\infty)\big)^2$ and $\p_1+\p_2=\nu(X)$. 

Consider the (right-continuous) ``distribution function'' $F$ of the function $f$ with respect to the measure $\nu$,  
%and the generalized inverse $F^{-1}$ 
defined by the formula%s 
\begin{equation*}
	F(t):=\nu\big(f^{-1}([0,t])\big)=\nu\big(\{x\in X\colon f(x)\le t\}\big)
\end{equation*}
for $t\in(-\infty,\infty]$, and let 
\begin{equation}\label{eq:q,t}
	%F^{-1}(u):=\inf F^{-1}([u,\infty))=
	s:=\inf\{t\in[0,\infty]\colon F(t)\ge \p_1\}\in[0,\infty],\quad\text{so that}\quad F(s-)\le \p_1\le F(s). 
\end{equation}
%for $u\in[0,\nu(X)]$. %; recall here that $\inf\emptyset=\infty$. 
%Next, 
Next, let $D:=f^{-1}(\{s\})$, and then  
let $D_1$ be any set in $\Si$ such that 
$D_1\subseteq D$ and $\nu(D_1)=\p_1-F(s-)$; 
such a set $D_1$ exists by \eqref{eq:p}, in view of the inequalities in \eqref{eq:q,t} and the equality $\nu(D)=F(s)-F(s-)$. 
Finally, let $B_1:=f^{-1}([0,s))\cup D_1$ and $B_2:=X\setminus B_1=f^{-1}((s,\infty))\cup(D\setminus D_1)$. Then, obviously, $(B_1,B_2)\in\PP_2$. Also, because $D_1\subseteq D=f^{-1}(\{s\})$, the sets $f^{-1}([0,s))$ and $D_1$ are disjoint and hence $\nu(B_1)=\nu\big(f^{-1}([0,s))\big)+\nu(D_1)=F(s-)+[\p_1-F(s-)]=\p_1$, so that $\nu(B_2)=\nu(X)-\nu(B_1)=\p_2$. Therefore, $(B_1,B_2)\in\PP_{\nu,\pp}$. 
Moreover, $B_1\subseteq f^{-1}([0,s])$ and $B_2\subseteq f^{-1}([s,\infty))$, so that $\sup_{B_1}f\le s\le\inf_{B_2}f$, and thus \eqref{eq:inf>sup} holds, for $k=2$. 
%
%%Take any $i\in\{0,\dots,k\}$. Let 
%%\begin{equation}\label{eq:q,t}
%%	\q_i:=\sum_{j=1}^i \p_j\quad\text{and}\quad t_i:=F^{-1}(\q_i),\quad\text{so that}\quad F(t_i-)\le \q_i\le F(t_i),  
%%\end{equation}
%%$\q_i\in G$, 
%%$\q_0=0$, and 
%%\begin{equation}\label{eq:t_i}
%%	0=t_0\le t_1\le\cdots\le t_k\le\infty. 
%%\end{equation}
%Next, let $D_{i,-}$ be any set in $\Si$ such that 
%\begin{equation}\label{eq:i+}
%	D_{i,-}\subseteq D_i:=f^{-1}(\{t_i\})\quad\text{and}\quad \nu(D_{i,-})=\q_i-F(t_i-); 
%\end{equation}
%such a set $D_{i,-}$ exists by \eqref{eq:p}, in view of the inequalities in \eqref{eq:q,t} and the equality $\nu(D_i)=F(t_i)-F(t_i-)$. 
%% 
%Next, let   
%\begin{equation}
%	D_{i,+}:=D_i\setminus D_{i,-}; 
%\end{equation}
%then 
%\begin{equation}\label{eq:i-}
%	\nu(D_{i,+})=F(t_i)-\q_i. 
%\end{equation}
%
%Now take any $i\in\{1,\dots,k\}$. Let 
%\begin{equation}
%D_{i-1,i}:=f^{-1}\big((t_{i-1},t_i)\big)\quad\text{and}\quad B_i:=D_{i-1,+}\cup D_{i-1,i}\cup D_{i,-}, 
%\end{equation}
%so that 
%\begin{equation}\label{eq:i-1,i}
%	\nu(D_{i-1,i})=F(t_i-)-F(t_{i-1}).  
%\end{equation}
%Note that the sets $D_{i-1,+}, D_{i-1,i},D_{i,-}$ are pairwise disjoint. 
%
%So, in view of \eqref{eq:i-}, \eqref{eq:i-1,i}, the second part of \eqref{eq:i+}, and the definition of $\q_i$ in \eqref{eq:q,t}, we conclude that $\nu(B_i)=\p_i$, for all $i\in\{1,\dots,k\}$; that is, by \eqref{eq:P_nu}, $Q:=(B_1,\dots,B_k)\in\PP_{\nu,\pp}$. 
%
%Moreover, for each $i\in\{1,\dots,k\}$ we have %$B_i\ne\emptyset$ (since $\nu(B_i)=\p_i>0$) and 
%$f(B_i)\subseteq[t_{i-1},t_i]$. So, \eqref{eq:inf>sup} follows by \eqref{eq:t_i}. 
This completes the proof of Proposition~\ref{prop:decr}. 
\end{proof}

\begin{proof}\emph{of Theorem~\ref{th:opt score}.\ \ } 
Let $\Pi:=\bigcup_{J\subseteq[k]}\Pi_J$, where $\Pi_J$ stands for the set of all permutations of the set $J$. Let $\T:=\Pi\times[k]\times[k]$. 
%
%Let $\T$ denote the set of all triples $(i,j,\pi)$ such that $(i,j)\in[k]^2$ and $\pi\in\Pi_J$ for some $J\subseteq[k]$, where $\Pi_J$ stands for the set of all permutations of the set $J$. 

Take any partition $Q$ as in Proposition~\ref{prop:decr} and any partition $P=(A_1,\dots,A_k)\break
\in\PP_{\nu,\pp}$. Introduce $C_{i,j}:=A_i\cap B_j$ for $(i,j)\in[k]^2$. A crucial observation is that the matrix $\big(\nu(C_{i,j})\big)_{i,j\in[k]}$ is nonnegative and sum-symmetric, and so, by Theorem~\ref{th:sum-symm}, 
\begin{equation*}
	\nu(C_{i,j})=\sum_{J\subseteq[k]}\sum_{\pi\in\Pi_J}w_\pi \ii{j=\pi(i)\in J} 
\end{equation*}
for all $(i,j)\in[k]^2$, 
where %$\Pi_J$ denotes the set of all permutations of the set $J$ and 
the $w_\pi$'s are some %nonnegative real 
numbers in $G\cap[0,\infty)$. %such that $\sum_{\pi\in\Pi}w_\pi=1$. 
Therefore and in view of \eqref{eq:p}, 
for each $(i,j)\in[k]^2$ there is a partition $(C_{\pi;i,j})_{\pi\in\Pi}$ of the set $C_{i,j}$ such that 
for each triple $(\pi,i,j)\in\T$ we have 
%such that $(i,j)\in[k]^2$ and $\pi\in\Pi_J$ for some $J\subseteq[k]$, 
$C_{\pi;i,j}\in\Si$ and
\begin{equation}\label{eq:nu()}
\nu(C_{\pi;i,j})=w_\pi \ii{j=\pi(i)\in J},  
\end{equation}
where, for any given $\pi\in\Pi$, the set $J\subseteq[k]$ is uniquely determined by the condition $\pi\in\Pi_J$. 
%for some $J\subseteq[k]$.  
% and $(C_{\pi;i,j})_{\pi\in\Pi}$ is a partition of $C_{i,j}$ for each $(i,j)\in[k]^2$, 
Hence, 
\begin{equation}\label{eq:nu,mu}
	\nu(C_{i,j})=\sum_{\pi\in\Pi}\nu(C_{\pi;i,j}) \quad\text{and}\quad \mu(C_{i,j})=\sum_{\pi\in\Pi}\mu(C_{\pi;i,j}). 
\end{equation} 
For each triple $(\pi,i,j)\in\T$, let  
\begin{equation*}\label{eq:r}
	r_{\pi;i,j}:=
	\left\{
	\begin{alignedat}{2}
	&\frac{\mu(C_{\pi;i,j})}{\nu(C_{\pi;i,j})}=\frac1{\nu(C_{\pi;i,j})}\int_{C_{\pi;i,j}}f\,d\nu
	\quad&&\text{if}\ \; \nu(C_{\pi;i,j})\ne0, \\ 
	&\sup_{B_j}\,f
	\quad&&\text{otherwise},  
	\end{alignedat}
	\right.
\end{equation*}
so that 
\begin{equation}\label{eq:mu=r nu}
\mu(C_{\pi;i,j})=r_{\pi;i,j}\,\nu(C_{\pi;i,j}); 	
\end{equation}
also, in view of the set inclusions $C_{\pi;i,j}\subseteq C_{i,j}\subseteq B_j$, we have $\inf_{B_j}\,f\le r_{\pi;i,j}\le\sup_{B_j}\,f$. 
%if $|C_{\pi;i,j}|\ne0$; otherwise, define $r_{\pi;i,j}$ arbitrarily, but so that one have 

Therefore, in view of inequalities \eqref{eq:inf>sup}, we now arrive at the second important point in this proof: that 
%, in view of the set inclusions $C_{\pi;i,j}\subseteq C_{i,j}\subseteq B_j$ for $(\pi;i,j)\in\T$ and inequalities \eqref{eq:inf>sup}, 
for all triples $(\pi,i_1,j_1)$ and $(\pi,i_2,j_2)$ in $\T$ we have the implication 
\begin{equation}\label{eq:r incr}
	j_1<j_2\implies r_{\pi;i_1,j_1}\le r_{\pi;i_2,j_2}.    
\end{equation}

By \eqref{eq:s}, \eqref{eq:nu,mu}, \eqref{eq:mu=r nu}, %{eq:r}, 
and \eqref{eq:nu()}, 
\begin{equation*}\label{eq:s(P)}
\begin{aligned}
	s(P)=\sum_{i\in[k]} s_i\mu(A_i)
	&=\sum_{i,j\in[k]} s_i\mu(C_{i,j}) \\ 
	&=\sum_{J\subseteq[k]}\sum_{\pi\in\Pi_J}\sum_{i,j\in[k]} s_i\mu(C_{\pi;i,j}) \\ 
	&=\sum_{J\subseteq[k]}\sum_{\pi\in\Pi_J}\sum_{i,j\in[k]} s_i r_{\pi;i,j}\nu(C_{\pi;i,j}) \\ 
	&=\sum_{J\subseteq[k]}\sum_{\pi\in\Pi_J}\sum_{i,j\in[k]} s_i r_{\pi;i,j}w_\pi \ii{j=\pi(i)\in J} \\ 
%&	=\sum_{J\subseteq[k]}\sum_{\pi\in\Pi_J}w_\pi\sum_{i\in J} s_i r_{\pi;i,\pi(i)} \\ 
&	=\sum_{J\subseteq[k]}\sum_{\pi\in\Pi_J}w_\pi\sum_{j\in J} s_{\pi^{-1}(j)}  r_{\pi;\pi^{-1}(j),j}.  
\end{aligned}	
\end{equation*}
%by (1), where 
%\begin{equation*}
%	a_{\pi;j}:=r_{\pi;\pi^{-1}(j),j}. \tag{4}
%\end{equation*}
Similarly to %\eqref{eq:s(P)}
this, %with $Q=(B_1,B_2,B_3)$ as before,  
we have 
\begin{align*}
	s(Q)=\sum_{i\in[k]} s_j\mu(B_j)
	&=\sum_{i,j\in[k]} s_j\mu(C_{i,j}) \\ 
%	&=\sum_{J\subseteq[k]}\sum_{\pi\in\Pi_J}\sum_{i,j\in[k]} s_j\mu(C_{\pi;i,j}) \\ 
%	&=\sum_{J\subseteq[k]}\sum_{\pi\in\Pi_J}\sum_{i,j\in[k]} s_i r_{\pi;i,j}\nu(C_{\pi;i,j}) \\ 
%	&=\sum_{J\subseteq[k]}\sum_{\pi\in\Pi_J}\sum_{i,j\in[k]} s_i r_{\pi;i,j}w_\pi \ii{j=\pi(i)\in J} \\ 
%&	=\sum_{J\subseteq[k]}\sum_{\pi\in\Pi_J}w_\pi\sum_{i\in J} s_j r_{\pi;i,\pi(i)} \\ 
&	=\sum_{J\subseteq[k]}\sum_{\pi\in\Pi_J}w_\pi\sum_{j\in J} s_j  r_{\pi;\pi^{-1}(j),j},   
\end{align*}
with the only difference that $s_i$ in $\sum_{i,j\in[k]} s_i\mu(C_{i,j})$ and in the two subsequent expressions in multi-line display \eqref{eq:s(P)} is now replaced by $s_j$. 

So, to compete the proof of Theorem~\ref{th:opt score}, it suffices to show that 
\begin{equation}\label{eq:s u}
	\sum_{j\in J} s_j u_j\ge\sum_{j\in J} s_{\si(j)} u_j
\end{equation}
for any permutation $\si\in\Pi_J$, where $u_j:=r_{\pi;\pi^{-1}(j),j}$. 
Since any permutation can be obtained from the identity permutation by finitely many inversions, it is enough to verify \eqref{eq:s u} in the case when the cardinality of $J$ is $1$ or $2$, so that $J=\{j,m\}$ for some $j,m$ in $[k]$. Then  
%is the product of transpositions, wlog the permutation $\si$ in \eqref{eq:s u} is a transposition of two elements, say $j$ and $m$, of $J$, and in this case 
\eqref{eq:s u} can be rewritten as $s_ju_j+s_mu_m\ge s_mu_j+s_ju_m$ or, equivalently, as $(s_j-s_m)(u_j-u_m)\ge0$, which is true -- because, by 
\eqref{eq:s incr} and \eqref{eq:r incr}, $s_j$ and $u_j=r_{\pi;\pi^{-1}(j),j}$ are each nondecreasing in $j\in J$. 
%
% Following the lines of the above proof, we can see that the optimal partition is unique, up to sets of zero Lebesgue measure. 
%
This concludes the proof of Theorem~\ref{th:opt score}. 
\end{proof}

%\newpage

\section{Optimal resource allocation} 
Theorem~\ref{th:opt score}, appropriately interpreted, provides a solution to an 
%the following 
optimal resource allocation (ORA) 
problem. 
For simplicity, let us state here this problem and its solution for the ``discrete'' setting, corresponding to alternative \eqref{II} on page~\pageref{I}. 
The ORA problem is as follows. 
%Let us interpret the finite set $X$ as a population of individuals. 

\begin{itemize}
	\item Each member $x$ of a finite set $X$ is to be subjected to exactly one of $k$ treatments, labeled by $1,\dots,k$, with potencies $s_1,\dots,s_k$ and available in quantities $q_1,\dots,q_k$, respectively. 
	\item In accordance with condition \eqref{eq:s incr}, we assume that the potencies $s_1,\dots,s_k$ are real numbers such that $s_1\le\dots\le s_k$; that is, the $k$ treatments are enumerated according to their potencies, from the lowest to the highest.  Potencies are allowed to take negative values, corresponding to negative treatment effects.
	\item 
%	In accordance with alternative \eqref{II} on page~\pageref{I} %(say) 
%	and condition  \eqref{eq:sum q}, let us also assume that the 
	In this ``discrete'' setting, the available quantities $q_1,\dots,q_k$ of treatments $1,\dots,k$ are nonnegative integers such that the total of the quantities $q_1,\dots,q_k$ equals the number $\nu(X)$ of the members of the set $X$. 
	\item For each member $x$ of the set $X$, the effect of any treatment $i\in[k]$ is proportional to the potency $s_i$ of the treatment, with a proportionality coefficient $f(x)\in[0,\infty)$, so that the just mentioned effect is $f(x)s_i$. It is then natural to refer to $f(x)$ as the responsiveness of member $x$ to treatment. 
	\item For each $i\in[k]$, let $A_i$ denote the set of all members $x$ of the set $X$ assigned to treatment $i$, so that $P:=(A_1,\dots,A_k)$ is a partition of $X$. This partition represents a \emph{treatment allocation}. In accordance with what has been said, we only consider ``feasible'' treatment allocations, that is, the ones satisfying the conditions $\nu(A_i)=\p_i$ for all $i\in[k]$; cf.\ \eqref{eq:P_nu} (recall that here $\nu$ stands for the counting measure). 
Letting now 
\begin{equation*}
\mu(A):=\int_A f\,d\nu=\sum_{x\in A}f(x) 	
\end{equation*}
for any set $A\subseteq X$, we see that the overall effect of a treatment allocation $P=(A_1,\dots,A_k)$ will then be  
\begin{equation*}
	\sum_{i\in[k]}\sum_{x\in A_i}f(x)s_i
	=\sum_{i\in[k]}s_i\mu(A_i)=s(P), 
\end{equation*}
in accordance with \eqref{eq:s}. 
\end{itemize}
 
Now Theorem~\ref{th:opt score} tells us that the overall effect $s(P)$ of a treatment allocation $P=(A_1,\dots,A_k)$ will be the largest possible if members of the set $X$ with higher responsiveness are assigned to higher-potency treatments. More specifically, for the optimal treatment allocation, $q_k$ members $x$ of the set $X$ with the highest values of responsiveness $f(x)$ are selected to constitute the set $A_k$ and thus to receive treatment $k$, of the highest-potency, $s_k$; then $q_{k-1}$ members of the remaining set $X\setminus A_k$ with the highest values of responsiveness are selected to constitute the set $A_{k-1}$ and thus to receive treatment $k-1$, of the second highest-potency, $s_{k-1}$; etc. 

While this %conclusion %may seem 
solution to this ORA 
%optimal resource allocation %(ORA) 
problem 
appears to agree with intuition, we saw that it takes some effort to prove it rigorously, by using the decomposition of nonnegative sum-symmetric matrices provided by Theorem~\ref{th:sum-symm}.

Let us now provide a few possible specific interpretations of the general ORA setting described above:  

\begin{enumerate}
	\item The set $X$ may be a human population to be vaccinated against a certain disease. Here, the treatments $1,\dots,k$ correspond to $k$ kinds of a vaccine, with potencies $s_1,\dots,s_k$. The total quantity of the available vaccine, $q_1+\dots+q_k$ units, is the same as the population size, so that each member of the population be able to receive exactly one unit of the vaccine. For each individual $x$ in the population, $f(x)$ is the individual's responsiveness to vaccination. The goal here is to maximize the overall vaccination effect $s(P)$. 
	\item Here $X$ is the set of workers of a certain specialty in an industrial company. 
Now the treatments $1,\dots,k$ correspond to $k$ kinds of equipment, with efficiencies $s_1,\dots,s_k$. The total quantity of the equipment units, $q_1+\dots+q_k$, is the same as the size of the set $X$ of workers, and each worker will be assigned to exactly one unit of the available equipment. For each worker $x$, $f(x)$ is the worker's individual productivity coefficient. The goal here is to maximize the overall production $s(P)$. 
	\item Now $X$ is a set of agricultural plots. The treatments $1,\dots,k$ correspond to $k$ grades of a fertilizer, with efficiencies $s_1,\dots,s_k$. The total quantity of the fertilizer units, $q_1+\dots+q_k$, is the same as the the number of plots, and each plot will receive exactly one unit of a fertilizer. For each plot $x$, $f(x)$ is the plot's responsiveness to fertilization. The goal here is to maximize the overall response $s(P)$ to the fertilization. 
		\item This is a ``non-atomic'' modification of the latter ``discrete'' scenario. Here $X$ is the set of points on an agricultural field, and the measure $\nu$ of a (measurable) part $A$ of $X$ is $c|A|$, where $c$ is a positive real number and $|A|$ is the area of $A$. The treatments $1,\dots,k$ again correspond to $k$ grades of a fertilizer, with efficiencies $s_1,\dots,s_k$. 
The field $X$ is partitioned into parts $A_1,\dots,A_k$ so that the part $A_j$ receive the $j$th grade of the fertilizer, for each $j=1,\dots,k$. 		
The corresponding quantities $q_1,\dots,q_k$ of the $k$ grades of the fertilizer may now take any nonnegative real values such that 
the total quantity of the fertilizer, $q_1+\dots+q_k$, equals $\nu(X)=c|X|$ so that the entire field be covered by %a layer of 
the fertilizer with the uniform density $c$ per unit area.  
For each point $x$ on the field, $f(x)$ is the corresponding local responsiveness to fertilization. The goal here is, again, to maximize the overall response $s(P)$ to the fertilization.
\end{enumerate}

In all these specific scenarios, the maximum overall effect occurs when higher levels of responsiveness are coupled with higher potencies, as specified in the general conclusion. 

\bigskip
\hrule
\bigskip

A search for articles containing the phrase ``optimal resource allocation'' in Google Scholar reveals about 35400 results. Optimal resource allocation (ORA) problems arise in a great variety of fields and a great variety of settings. A very small sample representing such problems includes ORA studies in 
biology~\cite{govern-wolde}, 
%in biology 
%compar-biology.pdf
%multi-cell.pdf, Optimal resource allocation in cellular sensing systems   PNAS.html 
computing~\cite{shahab-etal},  
%in computing 
%Optimal Resource Allocation in Clouds - IEEE Conference Publication.html 
%US6877035B2 - System for optimal resource allocation and planning for hosting computing services - Google Patents.html 
economics~\cite{arrow%,xianglan-jun
}, 
%
%Optimal allocation of resources in the material production sector and human capital sector - IEEE Conference Publication.html 
%Optimal Resource Allocation in an Imperfect Market Setting   Journal of Political Economy  Vol 69, No 6.html
electrical engineering \cite{seong-etal}, 
%seong-etal.pdf 
health care \cite{richter-etal}, 
%Applying Systems Engineering Principles in Improving Health Care Delivery   SpringerLink.html
%An Analysis of Optimal Resource Allocation for Prevention of Infection with Human Immunodeficiency Virus (HIV) in Injection Drug Users and Non-Users - Anke Richter, Margaret L. Brandeau, Douglas K. Owens, 1999.html 
information theory~\cite{li-goldsmith}, 
%in information theory 
%le-goldsmith.pdf
%tse-hanly.pdf
%wireless-netweorks.pdf 
operations research \cite{azaiez-bier}, 
%in operations research 
%Optimal resource allocation for security in reliability systems - ScienceDirect.html
%Optimal resource allocation across related channels - ScienceDirect.html
risk analysis \cite{bier-etal}, and 
%risk analysis in national defense 
%risk-in-defense.pdf
transportation \cite{dafermos-sparrow}. 
%transport-networks.pdf

Kantorovich was apparently the first to consider ORA problems systematically; see e.g.\ \cite[Section ``Linear programming'']{kantorovich} and \cite[page~240]{koopmans}. 
Methods used in the work by Kantorovich and his great many followers are analytical, based on separation of convex sets, with the feasible solutions being points in a finite- or infinite-dimensional linear space. 

On the other hand, the main tool used in the present paper is
the decomposition of nonnegative sum-symmetric matrices into nonnegative circuit
matrices, provided by Theorem~\ref{th:sum-symm}, whose %rather simple 
proof is rather combinatorial, and the feasible solutions in our setting are partitions, rather than points in linear spaces over $\R$. % of general form. 
It is hoped that the simple and rather general resource allocation model considered here, as well as the corresponding results, will be of use in a variety of specific applications. % settings/studies/situations. 

%Still, one may recall the similarity between the Birkhoff--von Neumann theorem (which 
%
%a firm mathematical foundation
%
%kantorovich-english.pdf

%
%
%[1]: https://en.wikipedia.org/wiki/Doubly_stochastic_matrix
%

% Preamble: \pgfplotsset{width=7cm,compat=1.12}

%\mathtoolsset{showonlyrefs,showmanualtags}

%% \appendix

%%  \bibliography{<your bibdatabase>}

%\begin{acknowledgements}
%If you'd like to thank anyone, place your comments here
%and remove the percent signs.
%\end{acknowledgements}

% BibTeX users please use one of
%\bibliographystyle{spbasic}      % basic style, author-year citations
\bibliographystyle{spmpsci}      % mathematics and physical sciences
%\bibliographystyle{spphys}       % APS-like style for physics
%\bibliography{}   % name your BibTeX data base
%
%\bibliographystyle{abbrv}
%%\bibliographystyle{amsplain}

\bibliography{P:/pCloudSync/mtu_pCloud_02-02-17/bib_files/citations10.13.18a}

\def\cprime{$'$} \def\polhk#1{\setbox0=\hbox{#1}{\ooalign{\hidewidth
  \lower1.5ex\hbox{`}\hidewidth\crcr\unhbox0}}}
  \def\polhk#1{\setbox0=\hbox{#1}{\ooalign{\hidewidth
  \lower1.5ex\hbox{`}\hidewidth\crcr\unhbox0}}}
  \def\polhk#1{\setbox0=\hbox{#1}{\ooalign{\hidewidth
  \lower1.5ex\hbox{`}\hidewidth\crcr\unhbox0}}} \def\cprime{$'$}
  \def\polhk#1{\setbox0=\hbox{#1}{\ooalign{\hidewidth
  \lower1.5ex\hbox{`}\hidewidth\crcr\unhbox0}}} \def\cprime{$'$}
  \def\polhk#1{\setbox0=\hbox{#1}{\ooalign{\hidewidth
  \lower1.5ex\hbox{`}\hidewidth\crcr\unhbox0}}} \def\cprime{$'$}
  \def\cprime{$'$}
\begin{thebibliography}{10}
\providecommand{\url}[1]{{#1}}
\providecommand{\urlprefix}{URL }
\expandafter\ifx\csname urlstyle\endcsname\relax
  \providecommand{\doi}[1]{DOI~\discretionary{}{}{}#1}\else
  \providecommand{\doi}{DOI~\discretionary{}{}{}\begingroup
  \urlstyle{rm}\Url}\fi

\bibitem{arrow}
Arrow, K.: Economic welfare and the allocation of resources for invention.
\newblock In: The Rate and Direction of Inventive Activity: Economic and Social
  Factors, pp. 609--626. National Bureau of Economic Research, Inc (1962)

\bibitem{azaiez-bier}
Azaiez, M.N., Bier, V.M.: Optimal resource allocation for security in
  reliability systems.
\newblock European J. Oper. Res. \textbf{181}(2), 773--786 (2007).
\newblock \doi{10.1016/j.ejor.2006.03.057}.
\newblock \urlprefix\url{https://doi.org/10.1016/j.ejor.2006.03.057}

\bibitem{bapat-ragh}
Bapat, R.B., Raghavan, T.E.S.: Nonnegative matrices and applications,
  \emph{Encyclopedia of Mathematics and its Applications}, vol.~64.
\newblock Cambridge University Press, Cambridge (1997).
\newblock \doi{10.1017/CBO9780511529979}.
\newblock \urlprefix\url{https://doi.org/10.1017/CBO9780511529979}

\bibitem{bier-etal}
Bier, V., Haphuriwat, N., Menoyo, J., Zimmerman, R., Culpen, A.: Optimal
  resource allocation for defense of targets based on differing measures of
  attractiveness.
\newblock Risk Anal. \textbf{28}, 763--770 (2008)

\bibitem{dafermos-sparrow}
Dafermos, S., Sparrow, F.T.: Optimal resource allocation and toll patterns in
  user-optimised transport networks.
\newblock Journal of Transport Economics and Policy \textbf{5}(2), 184--200
  (1971).
\newblock \urlprefix\url{http://www.jstor.org/stable/20052229}

\bibitem{dantzig85}
Dantzig, G.B., Eaves, B.C., Rothblum, U.G.: A decomposition and
  scaling-inequality for line-sum-symmetric nonnegative matrices.
\newblock SIAM J. Algebraic Discrete Methods \textbf{6}(2), 237--241 (1985).
\newblock \doi{10.1137/0606021}.
\newblock \urlprefix\url{https://doi.org/10.1137/0606021}

\bibitem{dudley-norv}
Dudley, R.M., Norvai\v{s}a, R.: Concrete functional calculus.
\newblock Springer Monographs in Mathematics. Springer, New York (2011).
\newblock \doi{10.1007/978-1-4419-6950-7}.
\newblock \urlprefix\url{https://doi.org/10.1007/978-1-4419-6950-7}

\bibitem{govern-wolde}
Govern, C.C., ten Wolde, P.R.: Optimal resource allocation in cellular sensing
  systems.
\newblock Proceedings of the National Academy of Sciences \textbf{111}(49),
  17,486--17,491 (2014).
\newblock \doi{10.1073/pnas.1411524111}.
\newblock \urlprefix\url{https://www.pnas.org/content/111/49/17486}

\bibitem{gravett}
Gravett, K.A.H.: Ordered abelian groups.
\newblock Quart. J. Math. Oxford Ser. (2) \textbf{7}, 57--63 (1956).
\newblock \doi{10.1093/qmath/7.1.57}.
\newblock \urlprefix\url{https://doi.org/10.1093/qmath/7.1.57}

\bibitem{hahn07}
Hahn, H.: \"{U}ber die nichtarchimedischen {G}r\"{o}{\ss}ensysteme.
\newblock Sitzungsberichte der {K}aiserlichen {A}kademie der {W}issenschaften,
  {W}ien, {M}athematisch -- {N}aturwissenschaftliche {K}lasse ({W}ien. {B}er.)
  \textbf{116}, 601--655 (1907)

\bibitem{hall67}
Hall Jr., M.: Combinatorial theory.
\newblock Blaisdell Publishing Co. Ginn and Co., Waltham, Mass.-Toronto,
  Ont.-London (1967)

\bibitem{kantorovich}
{Kantorovich}, L.V.: {My journey in science (proposed report to the Moscow
  Mathematical Society)}.
\newblock Russian Mathematical Surveys \textbf{42}, 233--270 (1987).
\newblock \doi{10.1070/RM1987v042n02ABEH001311}

\bibitem{koopmans}
Koopmans, T.C.: Concepts of optimality and their uses.
\newblock The American Economic Review \textbf{67}(3), 261--274 (1977).
\newblock \urlprefix\url{http://www.jstor.org/stable/1831399}

\bibitem{levi42}
Levi, F.W.: Ordered groups.
\newblock Proc. Indian Acad. Sci., Sect. A. \textbf{16}, 256--263 (1942)

\bibitem{li-goldsmith}
Li, L., {Goldsmith}, A.J.: Capacity and optimal resource allocation for fading
  broadcast channels---{P}art~ii: Outage capacity.
\newblock IEEE Transactions on Information Theory \textbf{47}(3), 1103--1127
  (2001).
\newblock \doi{10.1109/18.915667}

\bibitem{ann-comb}
Pinelis, I.: An extension of {H}all's theorem.
\newblock Ann. Comb. \textbf{6}(1), 103--106 (2002).
\newblock MR1923091

\bibitem{representant}
Pinelis, I.: A discrete mass transportation problem for infinitely many sites,
  and general representant systems for infinite families.
\newblock Math. Methods Oper. Res. \textbf{58}(1), 105--129 (2003).
\newblock MR2002566

\bibitem{richter-etal}
Richter, A., Brandeau, M., Owens, D.: An analysis of optimal resource
  allocation for prevention of infection with human immunodeficiency virus
  ({HIV}) in injection drug users and non-users.
\newblock Med. Decis. Making \textbf{19}, 167--179 (1999)

\bibitem{seong-etal}
{Seong}, K., {Mohseni}, M., {Cioffi}, J.M.: Optimal resource allocation for
  {OFDMA} downlink systems.
\newblock In: 2006 IEEE International Symposium on Information Theory, pp.
  1394--1398 (2006).
\newblock \doi{10.1109/ISIT.2006.262075}

\bibitem{shahab-etal}
Shahabuddin, J.S., {et al.}: System for optimal resource allocation and
  planning for hosting computing services (2005).
\newblock \urlprefix\url{https://patents.google.com/patent/US6877035}.
\newblock US Patent US6877035

\end{thebibliography}

%\bibliography{P:/mtu_pCloud_02-02-17/bib_files/citations12.13.12}

%\bibliography{C:/Users/ipinelis/Sync/mtu_Sync_01-19-17/bib_files/citations12.13.12}
%\bibliography{C:/Users/iosif/Sync/mtu_Sync_01-19-17/bib_files/citations12.13.12}
%\bibliography{C:/Users/Iosif/Dropbox/mtu/bib_files/citations12.13.12}

%%%% Non-BibTeX users please use
%%%\begin{thebibliography}{}
%%%%
%%%% and use \bibitem to create references. Consult the Instructions
%%%% for authors for reference list style.
%%%%
%%%\bibitem{RefJ}
%%%% Format for Journal Reference
%%%Author, Article title, Journal, Volume, page numbers (year)
%%%% Format for books
%%%\bibitem{RefB}
%%%Author, Book title, page numbers. Publisher, place (year)
%%%% etc
%%%\end{thebibliography}

\end{document}